\documentclass[11pt]{article}
\usepackage[margin=1in]{geometry}
\usepackage[T1]{fontenc}
\usepackage{lmodern,microtype}
\usepackage{amsmath,amssymb,amsthm,mathtools}
\usepackage{enumitem,needspace,booktabs,tabularx}
\usepackage[colorlinks=true,linkcolor=blue,citecolor=blue,urlcolor=blue]{hyperref}
\numberwithin{equation}{section}
\newtheorem{theorem}{Theorem}[section]
\newtheorem{proposition}[theorem]{Proposition}
\newtheorem{lemma}[theorem]{Lemma}
\newtheorem{corollary}[theorem]{Corollary}
\theoremstyle{definition}
\newtheorem{definition}[theorem]{Definition}
\theoremstyle{remark}
\newtheorem{remark}[theorem]{Remark}
\newcommand{\R}{\mathbb R}

\newcommand{\Hlog}{\mathcal H_{\log}}
\newcommand{\norm}[1]{\left\|#1\right\|}
\newcommand{\pair}[2]{\left\langle #1,#2\right\rangle}
\newcommand{\dH}[1]{\dot H^{#1}}
\newcommand{\Lop}{\mathcal L}
\DeclareMathOperator{\Rea}{Re}
\DeclareMathOperator{\Ima}{Im}
\title{Power-Type Nonlinear Schr\"odinger Equations\\with the Logarithmic Laplacian}
\author{Hichem Hajaiej\\[4pt]
\small Department of Mathematics\\
\small California State University, Los Angeles\\
\small Los Angeles, CA 90032, USA}
\date{}
\begin{document}
\maketitle
\vspace{-1.5em}

\begin{abstract}
We study the Cauchy problem for nonlinear Schr\"odinger equations with
Fourier dispersion $\log|\xi|^2$ and a pure-power nonlinearity. For
integer $m>N/2$ and sufficiently smooth powers, we establish local
well-posedness in $H^m$, mass conservation, strong-solution energy
conservation, and an amplitude-based continuation criterion. Annular
stationary phase and dyadic Strichartz estimates give global
small-data well-posedness and scattering for either sign in dimensions
$N\ge3$, including the focusing cubic equation. In the defocusing
case, a low-frequency estimate yields global weak solutions for
$2\le p\le3$ in the logarithmic form space intersected with
$L^{p+1}$. For focusing solutions, we prove failure of energy
coercivity and derive conditional lower bounds on finite-time Sobolev
blow-up. Energy conservation improves the mass-based lower bound by
introducing a stronger power and a logarithmic correction. The
occurrence of finite-time blow-up for large focusing data is not
established; the distinction between conditional rates and a
singularity construction is maintained throughout.
\end{abstract}


\section{Introduction and functional setting}

We consider the Cauchy problem
\begin{equation}\label{eq:main}
\begin{cases}
i\partial_tu=\Lop u+\mu |u|^{p-1}u,
 &(t,x)\in\R\times\R^N,\\
u(0,x)=u_0(x),&x\in\R^N,
\end{cases}
\qquad \Lop=\log(-\Delta),\quad p>1,
\end{equation}
where $u$ is complex-valued and $\mu\in\R$. The sign $\mu>0$ is
defocusing; the focusing equation is written with $\mu=-\kappa$, $\kappa>0$.
With the unitary Fourier transform,
\begin{equation}\label{eq:symbol}
\widehat{\Lop f}(\xi)=\log(|\xi|^2)\widehat f(\xi).
\end{equation}
The logarithmic Laplacian is connected to the fractional family through
the formal relation
\[
\log(-\Delta)=\left.\partial_s(-\Delta)^s\right|_{s=0}.
\]
Chen and Weth \cite{ChenWeth} developed its integral representation
and a functional framework for bounded-domain Dirichlet problems.
Their work provides background for the operator; the evolution results
proved here concern the whole-space nonlinear Cauchy problem
\eqref{eq:main}. Throughout this paper, the logarithm occurs in the
dispersion operator, rather than in a nonlinearity of the form
$u\log|u|^2$.

The logarithmic symbol is unbounded both above and below. At high
frequencies it grows more slowly than every positive power, while near
zero it has a negative singularity. Accordingly, the Sobolev, operator,
and logarithmic form domains must be distinguished. In particular,
high-frequency Sobolev regularity does not by itself remove the
low-frequency domain restriction for $\log(-\Delta)$. This issue
persists even though its Schr\"odinger group is unitary on every
Sobolev space.

For sufficiently regular nonlinearities, the semigroup contraction
method reduces local solvability to Sobolev product and composition
estimates. We follow that framework, as presented for dispersive
equations in \cite{Tao}, while keeping the low-frequency operator
domain separate from the space of mild solutions. At the global level,
unitarity alone is insufficient for the contraction argument. The
additional observation used here is that the logarithmic phase has a
nondegenerate Hessian on every fixed frequency annulus. This gives
annular time decay. The abstract theorem of Keel and Tao
\cite{KeelTao}, followed by spatial scaling and dyadic summation,
then supplies a global space-time bound with Sobolev regularity
$m>N/2$. That bound closes a small-data argument for both signs of
the nonlinearity.

Large-data energy methods have a different structure. For the
defocusing equation, the signed logarithmic kinetic form prevents an
immediate positive-energy estimate. We control its negative part
directly through the Fourier Duhamel formula. For $2\le p\le3$,
interpolation with the conserved mass reduces this control to a
linear Gronwall inequality. Smooth Hamiltonian approximations and
the compactness framework of Simon \cite{Simon} then yield global
weak solutions. This construction does not establish uniqueness in
the full energy class or global persistence of Sobolev regularity.

The focusing equation presents an additional difficulty: mass-preserving
concentration makes the energy unbounded below for every pure power
$p>1$. Fractional Schr\"odinger blow-up results, such as those of
Boulenger, Himmelsbach, and Lenzmann \cite{BHL} and Dinh
\cite{Dinh}, use positive-order dispersion in the range
$1/2<s<1$. They cannot be applied directly to
\eqref{eq:main}. We therefore distinguish the derivation of necessary
growth rates for an assumed finite lifespan from a proof that such
a lifespan actually occurs.

The results proved below are as follows.
\begin{enumerate}[label=(\roman*),leftmargin=2em,itemsep=2pt]
\item Local well-posedness in $H^m$, $m>N/2$, for smooth pure powers,
with mass conservation and an amplitude-based continuation criterion.
\item Frequency-localized dispersive estimates and global small-data
well-posedness and scattering, for either sign, in dimensions $N\ge3$.
\item Global defocusing weak solutions in the logarithmic form space
for $2\le p\le3$, without a Sobolev smoothness assumption on the power.
\item Focusing low-frequency estimates, failure of energy coercivity,
and dilation and variance identities under the requisite moment regularity.
\item Lower bounds on Sobolev growth conditional on a finite maximal
lifespan, including an energy-refined logarithmic correction.
\end{enumerate}
The local, global small-data, and weak-existence results are independent
statements with different assumptions; their exact hypotheses are
given in the corresponding theorems. The final two sections examine
focusing a priori estimates and conditional singularity rates.
No theorem asserting the occurrence of finite-time blow-up is proved.
In particular, conditional blow-up bounds are not constructions of
blow-up solutions. The arguments make no claim of optimality or priority.

\subsection{Domains and solution concepts}

All function spaces are over $\R^N$ unless otherwise stated. Write
$|D|^s$ for the Fourier multiplier $|\xi|^s$. The self-adjoint realization
of $\Lop$ on $L^2$ has domain
\begin{equation}\label{eq:domain}
D(\Lop)=\left\{f\in L^2:
\int |\log(|\xi|^2)|^2|\widehat f(\xi)|^2\,d\xi<\infty\right\}.
\end{equation}
The logarithmic form space is
\begin{equation}\label{eq:formspace}
\Hlog=D(|\Lop|^{1/2}),\qquad
\norm{f}_{\Hlog}^2=
\norm{f}_2^2+\int |\log(|\xi|^2)||\widehat f(\xi)|^2\,d\xi.
\end{equation}
Neither inclusion $H^m\subset D(\Lop)$ nor $H^m\subset\Hlog$ follows
from high-frequency regularity alone; the issue is at $\xi=0$.
For example, the Fourier profile
\[
\widehat f(\xi)=\mathbf1_{\{0<|\xi|<e^{-1}\}}
|\xi|^{-N/2}\left(\log\frac1{|\xi|}\right)^{-1}
\]
defines an element of every $H^m$, but $K_-(f)=\infty$: the corresponding
radial integrals are respectively $\int_1^\infty s^{-2}\,ds$ and
a constant multiple of $\int_1^\infty s^{-1}\,ds$.

Define
\begin{align}
K_+(f)&=\int_{|\xi|\ge1}\log(|\xi|^2)|\widehat f(\xi)|^2\,d\xi,\label{eq:Kplus}\\
K_-(f)&=\int_{|\xi|<1}-\log(|\xi|^2)|\widehat f(\xi)|^2\,d\xi.\label{eq:Kminus}
\end{align}
For $f\in\Hlog$, the signed kinetic form is $K(f)=K_+(f)-K_-(f)$.
The formal conserved quantities are
\begin{equation}\label{eq:energy}
M(f)=\norm{f}_2^2,\qquad
E_\mu(f)=\frac12K(f)+\frac{\mu}{p+1}\norm{f}_{p+1}^{p+1}.
\end{equation}

Set $U(t)=e^{-it\Lop}$, so
\begin{equation}\label{eq:group}
\widehat{U(t)f}(\xi)=|\xi|^{-2it}\widehat f(\xi).
\end{equation}
The value of the multiplier at $\xi=0$ is irrelevant. Since it has
modulus one almost everywhere, $U(t)$ is a strongly continuous unitary
group on every $H^s$, $s\in\R$. It also preserves $D(\Lop)$ and
$\Hlog$ with their graph and form norms.

\begin{definition}
A mild $H^m$ solution on an interval $I$ containing zero is a function
$u\in C(I;H^m)$ satisfying
\begin{equation}\label{eq:duhamel}
u(t)=U(t)u_0-i\mu\int_0^tU(t-s)|u(s)|^{p-1}u(s)\,ds
\end{equation}
in $H^m$. A strong $L^2$ solution additionally belongs to
$C(I;D(\Lop))\cap C^1(I;L^2)$ and satisfies the equation in $L^2$.
\end{definition}

For weak solutions, put
\begin{equation}\label{eq:V}
V_p=\Hlog\cap L^{p+1},\qquad
\norm{f}_{V_p}=\norm{f}_{\Hlog}+\norm{f}_{p+1}.
\end{equation}
This is a reflexive Banach space, continuously and densely embedded in
$L^2$. The kinetic form
\begin{equation}\label{eq:kinform}
\mathfrak l(u,v)=\int\log(|\xi|^2)\widehat u(\xi)
\overline{\widehat v(\xi)}\,d\xi
\end{equation}
defines a continuous map $\Lop:\Hlog\to\Hlog^*$.
A weak solution satisfies \eqref{eq:main} in $V_p^*$, with this form
interpretation of $\Lop$, and has the initial trace $u_0$ in $L^2$.

\section{Local well-posedness and conservation laws}

Throughout the Sobolev theory, use the following sufficient smoothness
hypothesis:
\begin{equation}\label{eq:powerhyp}
m\in\mathbb N,\quad m>N/2,\qquad
p\in\{3,5,7,\ldots\}\ \text{or}\ p>m+1.
\end{equation}
The restriction on real noninteger powers is convenient and is not
claimed to be sharp.

\begin{lemma}[Nonlinear estimates]\label{lem:moser}
Under \eqref{eq:powerhyp}, $F(u)=|u|^{p-1}u$ is locally Lipschitz
from $H^m$ to itself. More precisely,
\begin{align}
\norm{F(u)}_{H^m}
&\le C\norm{u}_\infty^{p-1}\norm{u}_{H^m},\label{eq:tame}\\
\norm{|D|^mF(u)}_2
&\le C\norm{u}_\infty^{p-1}\norm{|D|^m u}_2,\label{eq:homtame}
\end{align}
and
\begin{align}\label{eq:tamediff}
\norm{F(u)-F(v)}_{H^m}\le C\bigl[&
(\norm{u}_\infty+\norm{v}_\infty)^{p-1}\norm{u-v}_{H^m}\notag\\
&+(\norm{u}_\infty+\norm{v}_\infty)^{p-2}
(\norm{u}_{H^m}+\norm{v}_{H^m})\norm{u-v}_\infty\bigr].
\end{align}
\end{lemma}
\begin{proof}
For odd $p=2k+1$, the nonlinearity is the polynomial
$u^{k+1}\overline u^{\,k}$. For real $p>m+1$, it is a
$C^{m+1}$ map on $\R^2$, and
$|D^jF(z)|\le C_j|z|^{p-j}$ for the required derivative orders.
The chain rule, H\"older's inequality, and the integer-order
Gagliardo--Nirenberg estimates give \eqref{eq:tame} and
\eqref{eq:homtame}. For example, a term with $k$ differentiated factors
whose orders sum to $m$ is bounded by
$C\norm{u}_\infty^{p-k}\norm{u}_\infty^{k-1}
\norm{|D|^m u}_2$.
For the difference, write
\[
F(u)-F(v)=\int_0^1DF(v+\theta(u-v))(u-v)\,d\theta
\]
and apply the tame product and composition estimates. These are the
standard Sobolev estimates used in the semigroup method; see
\cite[Chapter 3 and Appendix A]{Tao}.
\end{proof}

\begin{theorem}[Local theory]\label{thm:local}
Assume \eqref{eq:powerhyp}. For every $u_0\in H^m$ and $\mu\in\R$,
there is a unique maximal mild solution
\[
u\in C((-T_-,T_+);H^m).
\]
The solution map is locally Lipschitz on bounded sets over a common
sufficiently short interval. If $T_+<\infty$, then
\begin{equation}\label{eq:alternative}
\lim_{t\uparrow T_+}\norm{u(t)}_{H^m}=\infty.
\end{equation}
For $\mu\ne0$, finite $T_+$ also implies
\begin{equation}\label{eq:ampcriterion}
\int_0^{T_+}\norm{u(t)}_\infty^{p-1}\,dt=\infty.
\end{equation}
The corresponding statements hold backward in time.
\end{theorem}
\begin{proof}
Use $C([-T,T];H^m)$ with its supremum norm. Sobolev embedding and
Lemma~\ref{lem:moser}, together with unitarity of $U(t)$, give for the
Duhamel map $\Phi$
\[
\norm{\Phi(u)}_{C H^m}\le\norm{u_0}_{H^m}
+C|\mu|T R^p,
\quad
\norm{\Phi(u)-\Phi(v)}_{C H^m}
\le C|\mu|T R^{p-1}\norm{u-v}_{C H^m}
\]
on a ball of radius $R$. For $u_0\ne0$, choose
$R=2\norm{u_0}_{H^m}$ and then $T$ so that the map preserves the ball
and is a contraction. Zero data give the zero solution.

For two solutions bounded by $R$ on a compact interval, Gronwall's
inequality yields
\[
\norm{u(t)-v(t)}_{H^m}
\le e^{C_R|\mu||t|}\norm{u_0-v_0}_{H^m}.
\]
This gives uniqueness and dependence on the data. Restarting the local
construction at any time proves \eqref{eq:alternative}, since the
local lifespan is bounded below on bounded $H^m$ sets.

Finally, \eqref{eq:tame} gives
\begin{equation}\label{eq:ampgronwall}
\norm{u(t)}_{H^m}\le\norm{u_0}_{H^m}
\exp\left(C|\mu|\int_0^t\norm{u(s)}_\infty^{p-1}\,ds\right).
\end{equation}
If the integral in \eqref{eq:ampcriterion} were finite, this estimate
would contradict \eqref{eq:alternative}.
\end{proof}

\begin{proposition}[Mass and strong solutions]\label{prop:strong}
Every mild solution in Theorem~\ref{thm:local} conserves mass.
If also $u_0\in D(\Lop)$, then on its maximal interval
\[
u\in C(H^m\cap D(\Lop))\cap C^1(L^2),
\]
where $D(\Lop)$ has its graph norm, and $E_\mu(u(t))=E_\mu(u_0)$.
\end{proposition}
\begin{proof}
The interaction variable $v(t)=U(-t)u(t)$ satisfies
\[
v(t)=u_0-i\mu\int_0^tU(-s)F(u(s))\,ds,
\]
and is $C^1$ in $H^m$. Unitarity gives
\[
\frac{d}{dt}\norm{u(t)}_2^2
=2\Rea\left(-i\mu\int|u(t)|^{p+1}\right)=0.
\]

The powers in \eqref{eq:powerhyp} satisfy $p\ge2$. Thus
$F(u)\in H^m\cap L^1$ and
\begin{equation}\label{eq:Fone}
\norm{F(u)}_1=\norm{u}_p^p
\le\norm{u}_\infty^{p-2}\norm{u}_2^2.
\end{equation}
For $f\in H^m\cap L^1$,
\begin{equation}\label{eq:domainbound}
\norm{\Lop f}_2\le C_m\norm{f}_{H^m}+C_N\norm{f}_1.
\end{equation}
At high frequencies this follows from the logarithm's subpower growth;
at low frequencies use $\norm{\widehat f}_\infty\le C\norm{f}_1$ and
$\int_{|\xi|<1}|\log(|\xi|^2)|^2\,d\xi<\infty$.
The nonlinearity is continuous in $H^m\cap L^1$ along the solution.
The Duhamel formula in the graph norm therefore proves persistence of
$D(\Lop)$ and $C^1(L^2)$ regularity.

For this regularity, the derivative of the kinetic form is
$2\Rea\pair{\Lop u}{u_t}$. The potential energy is differentiable
along the $L^2$ path, since $u$ is locally bounded in $L^\infty$.
Consequently
\[
\frac{d}{dt}E_\mu(u)
=\Rea\pair{\Lop u+\mu F(u)}{u_t}
=\Rea\pair{i u_t}{u_t}=0.
\]
\end{proof}

\begin{remark}[The Bessel logarithm]
The local theorem, mass conservation, and amplitude criterion hold also
for $\log(I-\Delta)$. In that case $H^m\subset D(\log(I-\Delta))$
automatically, so the local solutions above are strong $L^2$ solutions
without an extra low-frequency assumption. The global dispersive
argument in the next section is specifically for $\log(-\Delta)$;
it is not being transferred to the Bessel logarithm.
\end{remark}

\section{Frequency-localized dispersion and small-data scattering}

\begin{lemma}[Annular dispersion]\label{lem:annular}
Let $\psi\in C_c^\infty(\R^N\setminus\{0\})$ be supported in a fixed
annulus. Then
\begin{equation}\label{eq:annular}
\norm{U(t)\psi(D)f}_\infty
\le C_\psi(1+|t|)^{-N/2}\norm{f}_1.
\end{equation}
\end{lemma}
\begin{proof}
The kernel is
\[
K_\psi(t,x)=(2\pi)^{-N}\int
e^{i(x\cdot\xi-t\log|\xi|^2)}\psi(\xi)\,d\xi.
\]
For $|t|\le1$ it is bounded by $C\norm{\psi}_1$. For larger times,
write the phase as $t\Phi_y(\xi)$ with $y=x/t$ and
$\Phi_y(\xi)=y\cdot\xi-\log|\xi|^2$. Its Hessian is the negative of
\[
D^2\log|\xi|^2=\frac2{|\xi|^2}I-
\frac4{|\xi|^4}\xi\otimes\xi.
\]
The radial eigenvalue is $-2/|\xi|^2$, while tangential eigenvalues are
$2/|\xi|^2$. They are uniformly nonzero on the annulus. The map
$\xi\mapsto2\xi/|\xi|^2$ is injective, so there is at most one
stationary point. Parameter-uniform stationary phase near it and
nonstationary integration by parts on its complement give
$\sup_x|K_\psi(t,x)|\le C|t|^{-N/2}$. The parameter $y$ need only be
considered in a fixed compact set; outside it the phase gradient is
uniformly separated from zero. This proves the estimate.
\end{proof}

Let $P_j$ be a homogeneous smooth dyadic partition, supported where
$|\xi|\sim2^j$.

\begin{proposition}[Global space-time bound]\label{prop:strichartz}
Let $N\ge3$, $m>N/2$, and $2\le q<\infty$. Then
\begin{equation}\label{eq:globalLinfty}
\norm{U(t)f}_{L^q(\R;L^\infty)}\le C\norm{f}_{H^m}.
\end{equation}
Moreover, for $G\in L^1(\R;H^m)$,
\begin{equation}\label{eq:forced}
\norm{\int_0^tU(t-s)G(s)\,ds}_{L^\infty_tH^m_x}
+\norm{\int_0^tU(t-s)G(s)\,ds}_{L^q_tL^\infty_x}
\le C\norm{G}_{L^1_tH^m_x}.
\end{equation}
\end{proposition}
\begin{proof}
Apply the abstract Strichartz theorem of Keel--Tao
\cite[Theorem 1.2]{KeelTao} to $U(t)\psi(D)$. Its $TT^*$ kernel has
cutoff $|\psi|^2$ and satisfies Lemma~\ref{lem:annular}; its $L^2$
energy bound is immediate. Hence the unit-annulus estimate holds for
$q,r\ge2$ with $2/q+N/r=N/2$, including the endpoint
$(q,r)=(2,2N/(N-2))$. For $N\ge3$ each finite $q\ge2$ has such a
finite $r$.

Rescaling only space, and using
$\log|2^j\xi|^2=2j\log2+\log|\xi|^2$, gives
\begin{equation}\label{eq:dyadicStr}
\norm{U(t)P_jf}_{L^q_tL^r_x}
\le C2^{2j/q}\norm{P_jf}_2.
\end{equation}
Bernstein's inequality therefore implies
\[
\norm{U(t)P_jf}_{L^q_tL^\infty_x}
\le C2^{jN/r+2j/q}\norm{P_jf}_2
=C2^{jN/2}\norm{P_jf}_2.
\]
The dyadic sum is controlled by
\begin{align*}
\sum_{j\le0}2^{jN/2}\norm{P_jf}_2&\le C\norm{f}_2,\\
\sum_{j>0}2^{jN/2}\norm{P_jf}_2&\le
\left(\sum_{j>0}2^{2j(N/2-m)}\right)^{1/2}
\left(\sum_{j>0}2^{2jm}\norm{P_jf}_2^2\right)^{1/2}
\le C_m\norm{f}_{H^m}.
\end{align*}
This proves \eqref{eq:globalLinfty}, first on a dense class and then by
completion. Thus the derivative loss in \eqref{eq:dyadicStr} is
absorbed by the stated Sobolev regularity.

For \eqref{eq:forced}, the $L^\infty H^m$ part follows from unitarity.
For the other part, Minkowski's inequality and
\eqref{eq:globalLinfty} give, on the positive half-line,
\[
\norm{\int_0^tU(t-s)G(s)\,ds}_{L^q_{t\ge0}L^\infty_x}
\le \int_0^\infty
\norm{\mathbf1_{\{t\ge s\}}U(t-s)G(s)}_{L^q_{t\ge0}L^\infty_x}\,ds
\le C\int_0^\infty\norm{G(s)}_{H^m}\,ds.
\]
The negative half-line is identical. No retarded endpoint estimate
beyond this $L^1$ forcing bound is needed.
\end{proof}

\begin{theorem}[Small-data global well-posedness and scattering]\label{thm:small}
Let $N\ge3$ and assume \eqref{eq:powerhyp}. For each fixed
$\mu\in\R$ there is $\varepsilon_*>0$ such that
$\norm{u_0}_{H^m}\le\varepsilon_*$ implies a unique global mild solution
\[
u\in C_b(\R;H^m)\cap L^{p-1}(\R;L^\infty),
\]
with
\begin{equation}\label{eq:smallbound}
\sup_{t\in\R}\norm{u(t)}_{H^m}
+\norm{u}_{L^{p-1}_tL^\infty_x}
\le C\norm{u_0}_{H^m}.
\end{equation}
The small-data solution map is Lipschitz in this space, and there are
$u_\pm\in H^m$ such that
\begin{equation}\label{eq:scattering}
\norm{u(t)-U(t)u_\pm}_{H^m}\longrightarrow0
\quad\text{as }t\to\pm\infty.
\end{equation}
The theorem applies to the focusing as well as the defocusing sign.
\end{theorem}
\begin{proof}
Set $q=p-1\ge2$ and
$\mathcal X=C_b(\R;H^m)\cap L^q(\R;L^\infty)$, with the sum norm.
Lemma~\ref{lem:moser} and H\"older's inequality yield
\begin{align}
\norm{F(u)}_{L^1H^m}&\le C\norm{u}_{\mathcal X}^{p},\label{eq:nonlinX}\\
\norm{F(u)-F(v)}_{L^1H^m}&\le
C(\norm{u}_{\mathcal X}+\norm{v}_{\mathcal X})^{p-1}
\norm{u-v}_{\mathcal X}.\label{eq:diffX}
\end{align}
For the second term of \eqref{eq:tamediff}, the time exponents are
$q/(p-2)$ and $q$; their reciprocals sum to one.

Proposition~\ref{prop:strichartz} gives
\[
\norm{\Phi(u)}_{\mathcal X}\le C_0\norm{u_0}_{H^m}
+C_1|\mu|\norm{u}_{\mathcal X}^p,
\]
and the analogous contraction estimate from \eqref{eq:diffX}.
On the ball of radius $2C_0\norm{u_0}_{H^m}$ these close when
$C|\mu|\norm{u_0}_{H^m}^{p-1}$ is sufficiently small. Banach's fixed
point theorem gives a solution on the entire time axis, its Lipschitz
dependence, and \eqref{eq:smallbound}. Local uniqueness from
Theorem~\ref{thm:local} gives uniqueness among all $C(\R;H^m)$ mild
solutions with these data.

By \eqref{eq:nonlinX}, $F(u)\in L^1(\R;H^m)$. Define
\[
u_+=u_0-i\mu\int_0^\infty U(-s)F(u(s))\,ds.
\]
This integral converges in $H^m$, and
\[
\norm{u(t)-U(t)u_+}_{H^m}
\le |\mu|\int_t^\infty\norm{F(u(s))}_{H^m}\,ds\longrightarrow0.
\]
The backward integral gives $u_-$.
\end{proof}

\begin{corollary}[The focusing cubic equation]
For $N\ge3$, $m>N/2$ an integer, and sufficiently small
$u_0\in H^m$, the focusing cubic equation has a unique global mild
solution and scatters in $H^m$. In particular,
\[
\sup_t\norm{u(t)}_{H^m}
+\left(\int_\R\norm{u(t)}_\infty^2\,dt\right)^{1/2}
\le C\norm{u_0}_{H^m}.
\]
\end{corollary}

\section{Global defocusing weak solutions}

The theorem in this section does not require \eqref{eq:powerhyp}.
The restriction $2\le p\le3$ arises from a low-frequency differential
inequality and is not claimed to be an optimal range.

\begin{theorem}[Global weak existence]\label{thm:weak}
Let $N\ge1$, $\mu>0$, and $2\le p\le3$. For every $u_0\in V_p$,
there is a global weak solution of \eqref{eq:main} satisfying
\begin{equation}\label{eq:weakclass}
u\in C(\R;L^2)\cap L^\infty_{\rm loc}(\R;V_p),\qquad
u_t\in L^\infty_{\rm loc}(\R;V_p^*).
\end{equation}
It conserves mass and obeys, for constants depending on the initial
mass, form norm, nonlinear energy, and the parameters,
\begin{equation}\label{eq:weakexp}
K_+(u(t))+K_-(u(t))+\norm{u(t)}_{p+1}^{p+1}
\le C_0e^{C_1|t|}
\end{equation}
for almost every $t$. No general uniqueness, persistence of $H^m$, or
energy equality is asserted for this weak solution class.
\end{theorem}

\subsection{The low-frequency estimate}

For a regular solution, write $B(t)=K_-(u(t))^{1/2}$ and
$X(t)=\norm{u(t)}_{p+1}^{p+1}$. Since the group preserves the
low-frequency weighted $L^2$ norm,
\begin{equation}\label{eq:Bduhamel}
B(t)\le B(0)+C_N|\mu|\int_0^t\norm{u(s)}_p^p\,ds
\quad(t\ge0).
\end{equation}
Here we used
\[
\norm{\mathbf1_{\{|\xi|<1\}}
\sqrt{-\log|\xi|^2}\,\widehat f}_2\le C_N\norm{f}_1,
\]
which follows from the Fourier $L^1$ bound and integrability of
$-\log|\xi|^2$ on the unit ball.

For $p\ge2$, interpolation gives
\begin{equation}\label{eq:interpp}
\norm{u}_p^p\le\norm{u}_2^{2/(p-1)}
\norm{u}_{p+1}^{(p+1)(p-2)/(p-1)}.
\end{equation}
In the defocusing case, conservation of energy yields
\begin{equation}\label{eq:defocoercive}
\frac12K_+(u(t))+\frac\mu{p+1}X(t)
=E_\mu(u_0)+\frac12B(t)^2,
\end{equation}
and hence $X+K_+\le C(1+B^2)$. Set
$\alpha=(p-2)/(p-1)$. For $2\le p\le3$, $0\le\alpha\le1/2$,
so \eqref{eq:Bduhamel}--\eqref{eq:defocoercive} reduce to
\[
B(t)\le B(0)+C\int_0^t(1+B(s)^2)^\alpha\,ds
\le B(0)+C\int_0^t(1+B(s))\,ds.
\]
Gronwall's inequality gives the exponential bounds needed for global
weak compactness. The same argument works backward in time.

\subsection{Approximation, compactness, and the initial trace}

\begin{proof}[Proof of Theorem~\ref{thm:weak}]
Choose a real even $\chi\in C_c^\infty$ with $0\le\chi\le1$ and
$\chi=1$ on the unit ball. Let $J_n=\chi(D/n)$, with support contained
in $|\xi|\le2n$. These operators are uniformly bounded on
$L^r$, $1<r<\infty$, are $L^2$ contractions, commute with $\Lop$,
and map $L^2$ into $H^k\cap L^\infty$ for every fixed $k$.

Choose $u_{0,n}\to u_0$ in $V_p$ with
$u_{0,n}\in D(\Lop)\cap\bigcap_kH^k$. Such an approximation is
obtained by smooth Fourier cutoffs away from zero and infinity.
Convergence in $\Hlog$ follows by dominated convergence. For the
$L^{p+1}$ component, the high cutoff is an approximate identity; the
low cutoff tends strongly to zero on every finite $L^r$, $1<r<\infty$.
The latter follows first for $L^1\cap L^r$ functions by scaling of the
convolution kernel, and then by density and uniform $L^r$ bounds.

Consider
\begin{equation}\label{eq:reg}
i\partial_tu_n=\Lop u_n+\mu J_nF(J_nu_n),\qquad u_n(0)=u_{0,n}.
\end{equation}
For fixed $n$, its nonlinearity is locally Lipschitz on $L^2$ and is
bounded on each $L^2$ ball. The interaction-variable argument gives
mass conservation because $J_n$ is self-adjoint:
\[
\pair{J_nF(J_nu_n)}{u_n}=\int|J_nu_n|^{p+1}\in\R.
\]
Thus the approximate solution is global in $L^2$.
The forcing is band-limited, belongs to $L^1$ for $p\ge2$, and hence
belongs to $D(\Lop)$ by the low-frequency bound in
\eqref{eq:domainbound}. This gives graph-domain persistence for each
fixed $n$. Its conserved regularized energy is
\begin{equation}\label{eq:regenergy}
E_n=\frac12K_+(u_n)-\frac12K_-(u_n)
+\frac\mu{p+1}\norm{J_nu_n}_{p+1}^{p+1}.
\end{equation}
The initial $E_n$, masses, and $K_-(u_{0,n})$ are uniformly bounded.

Put $B_n=K_-(u_n)^{1/2}$ and
$X_n=\norm{J_nu_n}_{p+1}^{p+1}$. The weighted Duhamel estimate uses
$\norm{\widehat{J_nF(J_nu_n)}}_\infty
\le C\norm{F(J_nu_n)}_1$. Interpolation
\eqref{eq:interpp}, with $\norm{J_nu_n}_2\le\norm{u_n}_2$, gives
\[
B_n(t)\le B_n(0)+C\int_0^t(1+X_n(s))^\alpha\,ds,
\qquad K_+(u_n)+X_n\le C(1+B_n^2).
\]
For a precise Gronwall argument, bound $B_n$ by its integral envelope
$Z_n$, which satisfies $Z_n'\le C(1+Z_n^2)^\alpha\le C(1+Z_n)$.
Consequently, on every finite interval $[-T,T]$,
\begin{equation}\label{eq:unbounds}
\sup_n\sup_{|t|\le T}
\bigl(\norm{u_n(t)}_{\Hlog}^2+X_n(t)\bigr)\le C_T,
\end{equation}
with exponential dependence on $T$.

The equations bound $(u_n)_t$ in $L^\infty(-T,T;V_p^*)$:
the kinetic term is bounded by the form norm and
\[
|\pair{J_nF(J_nu_n)}{\phi}|
\le\norm{J_nu_n}_{p+1}^p\norm{J_n\phi}_{p+1}
\le C_T\norm{\phi}_{V_p}.
\]
To obtain local strong compactness, note that the high-frequency
logarithmic space with norm
\[
\int\log(2+|\xi|^2)|\widehat f(\xi)|^2\,d\xi
\]
is compactly embedded in $L^2$ on sets of functions supported in a
fixed bounded region. Indeed, its Fourier tails are bounded by
$C/\log(2+R^2)$, giving uniform translation continuity; bounded
support gives spatial tightness. Multiplication by a fixed smooth
compact cutoff is bounded in this high-frequency logarithmic space.

After such localization, the time derivatives are bounded in $H^{-k}$
for $k>N/2$. For the low-frequency part of $\Lop u_n$, use
\[
\int_{|\xi|<1}|\log|\xi|^2|\,|\widehat u_n|\,d\xi
\le K_-(u_n)^{1/2}
\left(\int_{|\xi|<1}|\log|\xi|^2|\,d\xi\right)^{1/2};
\]
its inverse transform is bounded locally. The high-frequency part is
bounded in $H^{-k}$ by \eqref{eq:unbounds}. The nonlinear term has the
same negative-Sobolev bound, since $H^k\hookrightarrow L^{p+1}$.
The Aubin--Lions--Simon compactness theorem \cite{Simon} therefore gives
\begin{equation}\label{eq:localstrong}
u_n\to u\quad\text{in }L^2((-T,T)\times B_R)
\end{equation}
after extraction, for every finite $T,R$.

Also,
\begin{equation}\label{eq:smoothingerror}
\norm{(I-J_n)u_n}_{L^\infty(-T,T;L^2)}^2
\le\frac{C_T}{\log(n^2)}\longrightarrow0.
\end{equation}
Thus $J_nu_n\to u$ locally in $L^2$ and almost everywhere.
Its global $L^{p+1}$ bound identifies the weak limit
\[
F(J_nu_n)\rightharpoonup F(u)
\quad\text{in }L^{(p+1)/p}((-T,T)\times\R^N).
\]
Passing to the limit in \eqref{eq:reg}, using
$J_n\phi\to\phi$ in $L^{p+1}$ for test functions, gives the equation
in $V_p^*$. Weak lower semicontinuity gives $u\in L^\infty(-T,T;V_p)$
and \eqref{eq:weakexp}. Equicontinuity of pairings with $V_p$ tests
and convergence of $u_{0,n}$ identify the initial trace as $u_0$.

Finally, the Gelfand triple $V_p\subset L^2\subset V_p^*$ gives
$u\in C([-T,T];L^2)$ and its norm chain rule. Taking the real part,
\[
\frac{d}{dt}\norm{u}_2^2
=2\Rea\pair{-i\Lop u-i\mu F(u)}{u}=0,
\]
since $\mathfrak l(u,u)$ and $\int|u|^{p+1}$ are real.
A diagonal extraction over increasing $T$ gives the global solution.
\end{proof}

\begin{proposition}[Conditional uniqueness in the energy class]\label{prop:weakunique}
Two solutions in \eqref{eq:weakclass} with the same initial data agree
on any common interval on which both belong to
$L^{p-1}_tL^\infty_x$.
\end{proposition}
\begin{proof}
The pointwise estimate
$|F(z)-F(w)|\le C_p(|z|^{p-1}+|w|^{p-1})|z-w|$
and cancellation of the real kinetic form give
\[
\frac{d}{dt}\norm{u-v}_2^2
\le C_p|\mu|(\norm{u}_\infty^{p-1}+\norm{v}_\infty^{p-1})
\norm{u-v}_2^2.
\]
The chain rule can be taken in $V_p^*$; alternatively the assumed
amplitude integrability gives a mild $L^2$ formulation. Gronwall's
inequality finishes the proof. The extra space-time bound is not a
consequence of Theorem~\ref{thm:weak}.
\end{proof}

\Needspace{10\baselineskip}
\section{Focusing estimates and limitations of the energy method}

In this section $\mu=-\kappa<0$. For strong solutions from
Proposition~\ref{prop:strong}, write $E_0=E_{-\kappa}(u_0)$ and
$M=\norm{u_0}_2^2$.

\begin{proposition}[Low-frequency and kinetic bounds]\label{prop:focusbounds}
For $p\ge2$, whenever the indicated quantities are finite,
\begin{align}
K_-(u(t))^{1/2}
&\le K_-(u_0)^{1/2}
+C_N\kappa M\int_0^t\norm{u(s)}_\infty^{p-2}\,ds,\label{eq:focuslow}\\
K_+(u(t))
&\le2|E_0|+K_-(u(t))
+\frac{2\kappa M}{p+1}\norm{u(t)}_\infty^{p-1}.\label{eq:focushigh}
\end{align}
\end{proposition}
\begin{proof}
Combine \eqref{eq:Bduhamel} with \eqref{eq:Fone}. Rearrangement of
the focusing energy and
$\norm{u}_{p+1}^{p+1}\le M\norm{u}_\infty^{p-1}$ give the second
inequality. These estimates depend on amplitude information; they do
not provide it independently.
\end{proof}

\begin{proposition}[Failure of energy coercivity at every mass]\label{prop:noncoercive}
For every $M>0$, $p>1$, and $\kappa>0$, the focusing energy is
unbounded below on smooth functions of mass $M$, even among functions
with $K_-=0$. No estimate
\[
\norm{u}_{p+1}^{p+1}\le C(M)(1+K_+(u))
\]
holds on that class.
\end{proposition}
\begin{proof}
Choose a nonzero Schwartz function $\varphi$ with Fourier support in
$1<|\xi|<2$ and mass $M$. For $\lambda\ge1$, put
$\varphi_\lambda(x)=\lambda^{N/2}\varphi(\lambda x)$.
Then
\begin{gather*}
M(\varphi_\lambda)=M,\quad K_-(\varphi_\lambda)=0,\\
K_+(\varphi_\lambda)=K_+(\varphi)+2M\log\lambda,\qquad
\norm{\varphi_\lambda}_{p+1}^{p+1}
=\lambda^{N(p-1)/2}\norm{\varphi}_{p+1}^{p+1}.
\end{gather*}
Therefore
\[
E_{-\kappa}(\varphi_\lambda)
=\frac12K_+(\varphi)+M\log\lambda
-\frac\kappa{p+1}\lambda^{N(p-1)/2}\norm{\varphi}_{p+1}^{p+1}
\longrightarrow-\infty.
\]
Polynomial concentration dominates logarithmic kinetic growth, also
contradicting the proposed inequality.
\end{proof}

\begin{remark}
Proposition~\ref{prop:noncoercive} is not a blow-up theorem. It also
does not contradict Theorem~\ref{thm:small}, which assumes smallness
of an $H^m$ norm, not smallness of mass alone.
\end{remark}

\subsection{Dilation and variance identities}

The identities below are conditional on sufficient regularity and
moment control to differentiate the displayed pairings and pass the
commutator or integration-by-parts limits. They are not asserted for
every energy-class weak solution. In particular, no propagation of
these moment assumptions is used to claim a blow-up theorem.

Let
\begin{equation}\label{eq:dilation}
D(t)=\Ima\int\overline u\,x\cdot\nabla u\,dx.
\end{equation}
For the focusing equation, the exact identity is
\begin{equation}\label{eq:virial}
D'(t)=2M-\frac{\kappa N(p-1)}{p+1}\norm{u(t)}_{p+1}^{p+1}.
\end{equation}
Indeed, the linear commutator contributes
\[
\int \xi\cdot\nabla_\xi\log|\xi|^2\,|\widehat u|^2\,d\xi=2M.
\]
Writing $\rho=|u|^2$, the nonlinear contribution is
\[
\kappa\int\rho\,x\cdot\nabla(\rho^{(p-1)/2})\,dx
=\frac{\kappa(p-1)}{p+1}\int x\cdot\nabla(|u|^{p+1})\,dx
=-\frac{\kappa N(p-1)}{p+1}\norm{u}_{p+1}^{p+1}.
\]
Integration gives the space-time relation
\begin{equation}\label{eq:dilationintegral}
\frac{\kappa N(p-1)}{p+1}\int_0^T\norm{u(t)}_{p+1}^{p+1}\,dt
=2MT+D(0)-D(T).
\end{equation}

For $V(t)=\int|x|^2|u(t,x)|^2\,dx$, the variance identity, when
justified, is
\begin{equation}\label{eq:variance}
V'(t)=-4\Ima\int
\overline{\widehat u(t,\xi)}\,
\frac{\xi}{|\xi|^2}\cdot\nabla_\xi\widehat u(t,\xi)\,d\xi.
\end{equation}
This follows from Plancherel and
$\nabla_\xi\log|\xi|^2=2\xi/|\xi|^2$; the local nonlinearity makes
no contribution to $V'$. For comparison,
$D=-\Ima\int\overline{\widehat u}\,\xi\cdot\nabla_\xi\widehat u$.
Thus $V'=4D$, the relation used for the classical Laplacian, is absent.
If $xu$ and $|D|^{-1}u$ belong to $L^2$, the right side of
\eqref{eq:variance} is absolutely controlled by
\[
|V'(t)|\le4\norm{xu(t)}_2\norm{|D|^{-1}u(t)}_2.
\]

\section{Conditional lower bounds on finite-time blow-up}

Here a blow-up rate means a necessary lower bound on the Sobolev norm
\emph{assuming} the maximal lifespan is finite. No upper bound,
asymptotic equivalence, or singularity construction is inferred.

\begin{theorem}[Mass-based lower bound]\label{thm:basicrate}
Assume \eqref{eq:powerhyp}, $\mu\ne0$, and $u_0\in H^m\setminus\{0\}$.
If $T_+<\infty$, then, for $t<T_+$,
\begin{equation}\label{eq:basicrate}
\norm{|D|^m u(t)}_2
\ge c_{N,m,p}|\mu|^{-2m/[N(p-1)]}
\norm{u_0}_2^{-(2m-N)/N}
(T_+-t)^{-2m/[N(p-1)]}.
\end{equation}
\end{theorem}
\begin{proof}
Put $Y(t)=\norm{|D|^m u(t)}_2$, $A=\norm{u_0}_2>0$, and
$\theta_0=N/(2m)$. The interaction variable is $C^1(H^m)$, so $Y$
is locally absolutely continuous. Lemma~\ref{lem:moser} gives
\[
Y'\le C|\mu|\norm{u}_\infty^{p-1}Y.
\]
The Gagliardo--Nirenberg inequality and mass conservation imply
\[
\norm{u}_\infty\le C A^{1-\theta_0}Y^{\theta_0},
\quad
Y'\le K Y^{1+q_0},\quad
q_0=\theta_0(p-1),\quad
K=C|\mu|A^{(1-\theta_0)(p-1)}.
\]
Since mass is fixed, \eqref{eq:alternative} implies $Y(t)\to\infty$
as $t\uparrow T_+$. Integrating the inequality for $Y^{-q_0}$ gives
$Y(t)^{-q_0}\le q_0K(T_+-t)$, which is \eqref{eq:basicrate}.
\end{proof}

\Needspace{15\baselineskip}
\begin{theorem}[Energy-refined focusing lower bound]\label{thm:refinedrate}
Assume \eqref{eq:powerhyp}, $\mu=-\kappa<0$, and
$0\ne u_0\in H^m\cap D(\Lop)$. If $T_+<\infty$, then for all $t$
sufficiently close to $T_+$,
\begin{equation}\label{eq:refinedrate}
\boxed{\displaystyle
\norm{|D|^m u(t)}_2
\ge c(T_+-t)^{-\gamma}
\left[\log\frac e{T_+-t}\right]^{-\delta},}
\end{equation}
where
\begin{equation}\label{eq:rateexponents}
\gamma=\frac{m(p+1)}{N(p-1)}-\frac12,
\qquad \delta=\frac mN-\frac12.
\end{equation}
The constant and the closeness threshold depend on the parameters,
mass, and energy. This theorem presupposes, rather than proves,
$T_+<\infty$.
\end{theorem}
\begin{proof}
Let $Y=\norm{|D|^m u}_2$, $M=\norm{u_0}_2^2>0$, and
$X=\norm{u}_{p+1}^{p+1}$. Energy conservation and $K_-\ge0$ imply
\begin{equation}\label{eq:Xenergy}
X\le\frac{p+1}{2\kappa}(K_++2|E_0|).
\end{equation}
For $|\xi|\ge1$,
$\log|\xi|^2\le m^{-1}\log(1+|\xi|^{2m})$. Jensen's inequality
with probability measure $M^{-1}|\widehat u|^2\,d\xi$ gives
\begin{equation}\label{eq:logJensen}
K_+(u)\le\frac M m\log\left(1+\frac{Y^2}{M}\right).
\end{equation}
Thus $X\le C\log(e+Y)$.

Interpolation between $L^{p+1}$ and $\dH m$ gives
\[
\norm{u}_\infty\le C\norm{u}_{p+1}^{1-\theta}Y^\theta,
\qquad
\theta=\frac{N}{m(p+1)-N(p-1)/2}\in(0,1).
\]
Consequently
\[
\norm{u}_\infty^{p-1}\le
C[\log(e+Y)]^aY^b,
\qquad
a=\frac{(p-1)(1-\theta)}{p+1},\quad
b=(p-1)\theta>0.
\]
The same interaction-variable argument as in
Theorem~\ref{thm:basicrate} gives, almost everywhere,
\begin{equation}\label{eq:rateODE}
Y'\le C Y^{1+b}[\log(e+Y)]^a.
\end{equation}
Since $Y(t)\to\infty$, integration toward $T_+$ yields
\[
C(T_+-t)\ge
\int_{Y(t)}^\infty
\frac{d\eta}{\eta^{1+b}[\log(e+\eta)]^a}
\ge cY(t)^{-b}[\log(e+Y(t))]^{-a}.
\]
The last inequality follows already by integrating over $[Y,2Y]$ for
large $Y$. Inverting this increasing-power logarithmic relation gives
\[
Y(t)\ge c(T_+-t)^{-1/b}
\left[\log\frac e{T_+-t}\right]^{-a/b}.
\]
Direct computation shows
$1/b=m(p+1)/[N(p-1)]-1/2$ and $a/b=m/N-1/2$.
\end{proof}

\begin{corollary}[Three-dimensional cubic example]
For $N=3$, $p=3$, $m=2$, a focusing strong solution with a finite
maximal $H^2$ lifespan necessarily satisfies, near that time,
\[
\norm{\Delta u(t)}_2\ge
c(T_+-t)^{-5/6}
\left[\log\frac e{T_+-t}\right]^{-1/6}.
\]
Without the operator-domain/energy assumption, the mass-based bound
still gives $\norm{\Delta u(t)}_2\ge c(T_+-t)^{-2/3}$.
\end{corollary}

\Needspace{10\baselineskip}
\section{Scope, unresolved questions, and comparison with fractional NLS}

Theorems~\ref{thm:local} and \ref{thm:small} are uniqueness statements
in Sobolev classes. Theorem~\ref{thm:weak} is a global existence result
in a larger energy class; its construction does not establish
uniqueness there or global persistence of a local $H^m$ solution.
Proposition~\ref{prop:weakunique} states the additional space-time
condition under which energy-class uniqueness does follow.

For focusing data of arbitrary size, the paper does not prove global
Sobolev existence or finite-time blow-up. Negative energy alone is not
shown to force blow-up. The rate theorems assume a finite maximal
lifespan and cannot be used as evidence that such a lifespan occurs.

One tempting comparison is the family
\[
\Lop_s=\frac{(-\Delta)^s-I}{s},\qquad s>0,
\]
whose symbols converge pointwise to $\log|\xi|^2$ away from zero.
Its linear dilation contribution is
\[
\int\xi\cdot\nabla_\xi
\left(\frac{|\xi|^{2s}-1}{s}\right)|\widehat u|^2\,d\xi
=2\norm{|D|^s u}_2^2.
\]
On suitable fixed test functions this tends to $2M$, not a
positive-order derivative norm. The fractional blow-up theorems in
\cite{BHL,Dinh} assume $1/2<s<1$ and do not apply to logarithmic
dispersion. Passing to a formal operator limit therefore does not
transfer their blow-up conclusions.

The remaining occurrence problem is to find admissible data and a
mechanism proving $T_+<\infty$, for example a justified localized
functional with a closed finite-time forcing inequality, or a
constructed singularity profile with error control. Neither mechanism
is supplied here. This is a statement about the scope of the present
arguments, not a claim of nonexistence of blow-up or an exhaustive
classification of the literature.

\end{document}